\documentclass{amsart}

\newtheorem{theorem}{Theorem}[section]

\theoremstyle{definition}
\newtheorem{definition}[theorem]{Definition}

\theoremstyle{remark}
\newtheorem{remark}[theorem]{Remark}

\newtheorem{question}[theorem]{Question}

\usepackage[backrefs]{amsrefs}
\newtheorem{openproblem}[theorem]{Open Problem}
\DeclareMathOperator{\linspan}{span}
\usepackage{esint}

\usepackage{combelow}

\numberwithin{equation}{section}

\begin{document}

\title{New Directions in Harmonic Analysis on $L^1$}


\author{Daniel Spector}
\address{
Daniel Spector\hfill\break\indent
National Chiao Tung University\hfill\break\indent
Department of Applied Mathematics\hfill\break\indent
Hsinchu, Taiwan}
\address{Nonlinear Analysis Unit\hfill\break\indent 
Okinawa Institute of Science and Technology Graduate University\hfill\break\indent
1919-1 Tancha, Onna-son, Kunigami-gun\hfill\break\indent 
Okinawa, Japan}
\address{
National Center for Theoretical Sciences\hfill\break\indent 
National Taiwan University\hfill\break\indent
No. 1 Sec. 4 Roosevelt Rd.\hfill\break\indent 
Taipei, 106, Taiwan}
\address{
Washington University in St. Louis  \hfill\break\indent 
Department of Mathematics and Statistics\hfill\break\indent
One Brookings Drive\hfill\break\indent 
St. Louis, MO}
\curraddr{}
\email{dspector@math.nctu.edu.tw}
\thanks{}


\subjclass[2010]{Primary }

\date{}

\dedicatory{To Haim Brezis, on the occasion of his 75th birthday, with admiration and gratitude}

\begin{abstract}
The study of what we now call Sobolev inequalities has been studied for almost a century in various forms, while it has been eighty years since Sobolev's seminal mathematical contributions.  Yet there are still things we don't understand about the action of integral operators on functions.  This is no more apparent than in the $L^1$ setting, where only recently have optimal inequalities been obtained on the Lebesgue and Lorentz scale for scalar functions, while the full resolution of similar estimates for vector-valued functions is incomplete.  The purpose of this paper is to discuss how some often overlooked estimates for the classical Poisson equation give an entry into these questions, to the present state of the art of what is known, and to survey some open problems on the frontier of research in the area. 
\end{abstract}

\maketitle
\tableofcontents

\section{Introduction}\label{introduction}

One of the starting points for many interesting questions in harmonic analysis is a classical problem from the field of partial differential equations, the Poisson equation:  Given $f \in L^p(\mathbb{R}^d)$, find $u \in L^1_{loc}(\mathbb{R}^d)$ such that
\begin{align}
-\Delta u = f  \text{ in } \mathbb{R}^d \label{poisson}
\end{align}
in an appropriate sense.  Of course, to compute a distributional solution to \eqref{poisson} is not the difficult part of the problem, as we discuss\footnote{Our interest here is more than pedagogical, as our derivation gives rise to a new representation of the solution in the case $d=2$, which was originally proved by the author and R. Garg in \cite{GargSpector,GargSpector1}.}
 below in Section \ref{potentials}.  The main point is that with such a solution we want certain a priori estimates.  The three most basic estimates (and let us for the sake of discussion assume $d\geq 3$ here) one can ask for are the inequalities
\begin{align}
\|u\|_{L^q(\mathbb{R}^d)} &\leq C \|f\|_{L^p(\mathbb{R}^d)} ,\quad \quad \frac{1}{q} = \frac{1}{p}-\frac{2}{d}, \label{u} \\
\|\nabla u\|_{L^{q'}(\mathbb{R}^d)} &\leq C' \|f\|_{L^p(\mathbb{R}^d)}, \quad \quad \frac{1}{q'} = \frac{1}{p}-\frac{1}{d}, \text{ and }\label{gradu} \\
\|\nabla^2 u\|_{L^p(\mathbb{R}^d)} &\leq C'' \|f\|_{L^p(\mathbb{R}^d)}, \quad \quad \frac{1}{p} = \frac{1}{p}. \label{hessianu}
\end{align}
Then the classical results concerning \eqref{poisson} are that 
\begin{align*}
&\eqref{u} \quad \text{ holds whenever } \quad 1<p<d/2,\\
&\eqref{gradu} \quad \text{ holds whenever } \quad 1<p<d,\\
&\eqref{hessianu} \quad \text{ holds whenever } \quad 1<p<+\infty,
\end{align*}
while in the case $p=1$ one has a counterexample that shows all three of these inequalities are all false.  
This raises two natural questions. Firstly, one poses
\begin{question}\label{q1}
If one insists to obtain estimates for $u$, $\nabla u$, and $\nabla^2u$ in terms of $\|f\|_{L^1(\mathbb{R}^d)}$, what are the best possible spaces for such estimates?
\end{question}
We can give a fairly satisfying answer to this question by replacing the usual Lebesgue spaces with weak-type spaces, which are now commonplace in harmonic analysis.  In particular, observe that  Chebychev's inequality leads one to a natural quasi-norm of the functions we wish to estimate with the right scaling
\begin{align*}
\sup_{t>0} t |\{ |u|>t\}|^{(d-2)/d}  &\leq \|u\|_{L^{d/(d-2)}(\mathbb{R}^d)}, \\
\sup_{t>0} t |\{ |\nabla u|>t\}|^{(d-1)/d} &\leq \|\nabla u\|_{L^{d/(d-1)}(\mathbb{R}^d)}, \\
\sup_{t>0} t |\{ |\nabla^2 u|>t\}| &\leq \|\nabla^2 u\|_{L^{1}(\mathbb{R}^d)}.
\end{align*}
Thus, while it is not possible to control the right hand side of these inequalities by $\|f\|_{L^1(\mathbb{R}^d)}$, one might try to obtain estimates for the left hand side.  This can indeed be accomplished, that one has the inequalities
\begin{align}
\sup_{t>0} t |\{ |u|>t\}|^{(d-2)/d} &\leq C \|f\|_{L^1(\mathbb{R}^d)}  \label{u'}, \\
\sup_{t>0} t |\{ |\nabla u|>t\}|^{(d-1)/d}&\leq \tilde{C} \|f\|_{L^1(\mathbb{R}^d)}\label{gradu'},\\
\sup_{t>0} t |\{ |\nabla^2 u|>t\}| &\leq \tilde{\tilde{C}}  \|f\|_{L^1(\mathbb{R}^d)}.\label{hessianu'}
\end{align}
The estimates \eqref{u'} and \eqref{gradu'} can be obtained by some ideas in the 1956 paper of Zygmund \cite{Zygmund} detailing and extending some results of Marcinkiewicz, while the estimate \eqref{hessianu'} is in earlier work of Calder\'on and Zygmund from 1952 (see Lemma 2 in \cite{CalderonZygmund}).

This gives a fairly satisfying answer to Question \ref{q1}, and in fact it is the best one can hope for on the natural Lorentz\footnote{We discuss some of the value of the Lorentz spaces in Section \ref{Lorentz} below, where we also discuss stronger versions of the Lebesgue results presented here in the Introduction.} scale.  On the other hand, one could attempt to strengthen the hypothesis of the theorem for $p=1$ to obtain an analogous result to the case $p>1$, which can be asked as
\begin{question}\label{q2}
If one insists to obtain estimates on $u, \nabla u, \nabla^2 u$ in the Lebesgue spaces which scale correctly with $\|f \|_{L^1(\mathbb{R}^d)}$, what are the optimal assumptions to place on $f$ to ensure such estimates hold?
\end{question}
Concretely we are here asking what should one utilize for the right-hand-side in the inequalities
\begin{align*}
\|u\|_{L^{d/(d-2)}(\mathbb{R}^d)} &\leq \quad ?, \\
\|\nabla u\|_{L^{d/(d-1)}(\mathbb{R}^d)} &\leq \quad ?,\\
 \|\nabla^2 u\|_{L^1(\mathbb{R}^d)} &\leq \quad ?.
\end{align*}
Now for a replacement of \eqref{hessianu} it was subsequently understood that if one assumes $f \in \mathcal{H}^1(\mathbb{R}^d)$, the real Hardy space, then one has the estimate
\begin{align*}
 \|\nabla^2 u\|_{L^1(\mathbb{R}^d)} &\leq C''\|f\|_{\mathcal{H}^1(\mathbb{R}^d)},
\end{align*}
and for a number of reasons this is a satisfactory answer.  However, this also led to replacements of \eqref{u} and \eqref{gradu} in terms of the Hardy space,
\begin{align}
\|u\|_{L^{d/(d-2)}(\mathbb{R}^d)} &\leq C\|f\|_{\mathcal{H}^1(\mathbb{R}^d)} \label{hardyu} \\
\|\nabla u\|_{L^{d/(d-1)}(\mathbb{R}^d)} &\leq C'\|f\|_{\mathcal{H}^1(\mathbb{R}^d)},\label{hardygradu}
\end{align}
and the main new directions we are interested in here stem from the fact that these embeddings are not optimal.

To discuss this lack of optimality let us cast the problem in a slightly more general setting, introducing the Riesz potentials 
\begin{align*}
I_\alpha f(x) = \frac{1}{\gamma(\alpha)}\int_{\mathbb{R}^d} \frac{f(y)}{|x-y|^{d-\alpha}}\;dy,
\end{align*}
for $\gamma(\alpha)$ (defined in Section \ref{potentials}) such that they satisfy the semi-group property
\begin{align*}
I_\alpha I_\beta f = I_{\alpha+\beta} f
\end{align*}
for $\alpha, \beta \in (0,d)$ and $\alpha+\beta<d$ and $f$ sufficiently nice.  Then the solution to Poisson's equation in the case $d \geq 3$ is simply
\begin{align*}
u=I_2 f,
\end{align*}
while more generally we have a notion of integration in several dimensions which allows one to integrate a suitably decaying function by order $\alpha \in (0,d)$.  In contrast to iterated one dimensional integration, this integration is by construction radial, taking values of a function on spheres and giving them equal weight to the value of the potential at the center, closer spheres being weighted more heavily.

In this framework we can write one fundamental estimate from which one can easily obtain \eqref{u}, \eqref{gradu}, and \eqref{hessianu} in the case $p>1$, the following  theorem about integrals of the potential type due to S. Sobolev \cite{sobolev} in 1938 (see below in Section \ref{potentialestimates} for a deduction of the inequalities \eqref{u}, \eqref{gradu}, and \eqref{hessianu} from this theorem).
\begin{theorem} \label{Sobolev}
Let $0<\alpha<d$ and $1<p<d/\alpha$.  Then there exists a constant 
$C=C(p,\alpha,d)>0$ such that
\begin{align}
\|I_\alpha f\|_{L^{q}(\mathbb{R}^d)} \leq C 
\|f\|_{L^p(\mathbb{R}^d)} \label{sobolevineq}
\end{align} 
for all $f \in L^p (\mathbb{R}^d)$, where
\begin{align*}
\frac{1}{q} = \frac{1}{p} - \frac{\alpha}{d}.
\end{align*}
\end{theorem}
As in the case of the inequalities \eqref{u}, \eqref{gradu}, and \eqref{hessianu}, no such inequality can hold in the case $p=1$.  The counterexample here is the same as before, and as it is instructive for our discussion, let us here detail it.  Let us suppose one had such an inequality for all $f \in L^1(\mathbb{R}^d)$, and let $\{f_n\}$ be a sequence such that $f_n \overset{*}{\rightharpoonup} \delta_0$.  One can take, for example, $f_n(x) = \frac{1}{|B(0,\frac{1}{n})|} \chi_{B(0,\frac{1}{n})}(x)$.  Then
\begin{align*}
\|I_\alpha f_n\|_{L^{q}(\mathbb{R}^d)} \leq C \|f_n\|_{L^1(\mathbb{R}^d)} = C,
\end{align*}
while 
\begin{align*}
I_\alpha f_n \to I_\alpha \ast \delta_0 = \frac{1}{\gamma(\alpha)} \frac{1}{|x|^{d-\alpha}},
\end{align*}
as the Dirac delta is the identity for convolution (this can also be verified to hold almost everywhere by the Lebesgue differentiation theorem).  Thus Fatou's lemma would imply
\begin{align*}
\left\|  \frac{1}{\gamma(\alpha)} \frac{1}{|\cdot|^{d-\alpha}} \right\|_{L^q(\mathbb{R}^d)} \leq C,
\end{align*}
but as $q=d/(d-\alpha)$, this would mean
\begin{align*}
+\infty = \left(  \int_{\mathbb{R}^d} \frac{1}{|x|^d}\;dx \right)^{(d-\alpha)/d} =  \left(  \int_{\mathbb{R}^d} \left|\frac{1}{|x|^{d-\alpha}}\right|^{d/(d-\alpha)}\;dx \right)^{(d-\alpha)/d} \leq C\gamma(\alpha) ,
\end{align*}
which is absurd.

By making a stronger assumption, an extension of Theorem \ref{Sobolev} to the case $p=1$ was proved by E. Stein and G. Weiss in 1960 (see \cite{SteinWeiss}).  In particular their result implies the following
\begin{theorem}[Stein-Weiss]\label{steinweiss}
Let $\alpha \in (0,d)$.  There exists a constant 
$C=C(\alpha,d)>0$ such that
\begin{align}
\|I_\alpha f\|_{L^{d/(d-\alpha)}(\mathbb{R}^d)} \leq C 
\left(\|f\|_{L^1(\mathbb{R}^d)} + \|Rf\|_{L^1(\mathbb{R}^d;\mathbb{R}^d)}\right) \label{steinweisseqn}
\end{align}
for all $f \in L^1(\mathbb{R}^d)$ such that $Rf:=\nabla I_1 f \in L^1(\mathbb{R}^d;\mathbb{R}^d)$ 
\end{theorem}

\begin{remark}
The condition $f \in L^1(\mathbb{R}^d)$ and its vector-valued Riesz transform $Rf:=\nabla I_1 f \in L^1(\mathbb{R}^d;\mathbb{R}^d)$ is the original definition of the real Hardy space $\mathcal{H}^1(\mathbb{R}^d)$ in several variables due to E. Stein and G. Weiss in \cite{SteinWeiss}.  One now has a number of equivalent definitions, for example, in terms of maximal functions \cite{FeffermanStein} or via an atomic decomposition \cite{Coifman,Latter}.
\end{remark}

How can such a theorem hold in light of the failure of the inequality \eqref{sobolevineq} for $p=1$? Well, formally one has
\begin{align*}
R\delta_0 = c_d \frac{x\;\;\;\;}{|x|^{d+1}} \notin L^1_{loc}(\mathbb{R}^d),
\end{align*}
so that both the left-hand-side and the right-hand-side of the inequality blow up along such an approach sequence as constructed above.  

The fact that the $L^1$-norm of the Riesz transform term blows up, and not the $L^1$-norm of the function itself, might suggest a possible improvement to such an inequality in the removal of the term $\|f \|_{L^1(\mathbb{R}^d)}$.  However one also has a more convincing argument of this fact which comes from an inequality arising in PDEs.  In particular, let us recall that E. Gagliardo \cite{Gagliardo} and L. Nirenberg \cite{Nirenberg} had proved the inequality (for $d\geq 2$, the case $d=1$ being an easy consequence of the fundamental theorem of calculus):
\begin{align*}
\| u\|_{L^{d/(d-1)}(\mathbb{R}^d)} \leq C \|\nabla u \|_{L^1(\mathbb{R}^d;\mathbb{R}^d)}
\end{align*}
for all $u$ sufficiently nice.  The choice $u=I_1f$ yields
\begin{align*}
\|I_1 f\|_{L^{d/(d-1)}(\mathbb{R}^d)}  \leq C \|Rf\|_{L^1(\mathbb{R}^d;\mathbb{R}^d)},
\end{align*}
which improves \eqref{steinweisseqn}.  

Thus one can improve the estimate of Stein and Weiss for $I_1$, and so one wonders
\begin{question}
Can one make a similar improvement for $\alpha \in (0,d)$?
\end{question}

The first observation in this regard is that if $\alpha \in [1,d)$,  the semi-group property of the Riesz potentials, Sobolev's inequality \eqref{sobolevineq}, and the inequality of E. Gagliardo \cite{Gagliardo} and L. Nirenberg \cite{Nirenberg} imply
\begin{align*}
\|I_\alpha f\|_{L^{d/(d-\alpha)}(\mathbb{R}^d)}  &= \|I_{\alpha-1} I_1 f\|_{L^{d/(d-\alpha)}(\mathbb{R}^d)} \\
&\leq C \|I_1 f\|_{L^{d/(d-1)}(\mathbb{R}^d)} \\
&\leq C'\|Rf\|_{L^1(\mathbb{R}^d;\mathbb{R}^d)}.
\end{align*}
Thus one has such improvements for $\alpha \in [1,d)$.  Naturally this range of $\alpha$ contains the values $2$ and $1$, which are correspond to the amount of integration being performed in \eqref{hardyu} and \eqref{hardygradu}, respectively.  One immediately deduces an improvement to \eqref{hardyu}, while a similar improvement to \eqref{hardygradu} is a little more subtle.  In particular we require the boundedness of the Riesz transform 
\begin{align*}
R :L^p(\mathbb{R}^d) \to L^p(\mathbb{R}^d;\mathbb{R}^d) \text{ for } 1<p<+\infty,
\end{align*}
see e.g.~p.~33 in \cite{Sharmonic}.  This boundedness and the validity of the formula 
\begin{align*}
\nabla u(x) = \nabla I_1 (I_1f) \equiv R (I_1f)
\end{align*}
implies
\begin{align*}
\| \nabla u\|_{L^{d/(d-1)}(\mathbb{R}^d;\mathbb{R}^d)} \leq C \|I_1 f\|_{L^{d/(d-1)}(\mathbb{R}^d)},
\end{align*}
from which the result follows again from the case $\alpha=1$.

Thus we have seen that there are improvements to the inequalities \eqref{hardyu} and \eqref{hardygradu}, and to the potential mapping properties for any $\alpha \in [1,d)$.   This motivates
\begin{question}
Suppose $d \geq 2$.  Can one show that for $\alpha \in (0,1)$ there exists a constant $C=C(\alpha,d)>0$ such that
\begin{align}\label{potentialnodirac}
\|I_\alpha f \|_{L^{d/(d-\alpha)}(\mathbb{R}^d)} \leq C\|Rf\|_{L^1(\mathbb{R}^d;\mathbb{R}^d)}
\end{align}
for all $f \in C^\infty_c(\mathbb{R}^d)$ such that $Rf\in L^1(\mathbb{R}^d;\mathbb{R}^d)$?
\end{question}
\noindent
That we consider the case $d \geq 2$ here goes beyond the fact that the classical estimates for Poisson's equation we have recorded have no meaning, or would have to be suitably interpreted.  While the scaling would allow for such an inequality when $\alpha \in (0,1)$, in one dimension the estimate has a fundamental obstruction.  In particular, in this setting the Riesz transform collapses to the Hilbert transform, so that the inequality would read
\begin{align*}
\|I_\alpha f\|_{L^{1/(1-\alpha)}(\mathbb{R})} \leq C \| H f\|_{L^1(\mathbb{R})}.
\end{align*}
However, now the identity $H^2=-I$ and the boundedness of $H: L^p(\mathbb{R}) \to L^p(\mathbb{R})$ for $1<p<+\infty$ would imply
\begin{align*}
\|I_\alpha f\|_{L^{1/(1-\alpha)}(\mathbb{R})} \leq C \| f\|_{L^1(\mathbb{R})},
\end{align*}
which is precisely the estimate we have contradicted.

The estimate \eqref{potentialnodirac} extrapolates the overlooked estimates from Poisson's equation, that for the more general Riesz potentials of any order $\alpha \in (0,d)$, while one cannot show an $L^1$ estimate, perhaps one can show these $L^1$-type estimates. It is perhaps no surprise to the reader at this point that the estimate \eqref{potentialnodirac} is indeed valid.   As was discussed in the work of A. Schikorra, the author, and J. Van Schaftingen in \cite{SSVS}, one has a number of more classical approaches to prove the inequality, provided one knows the estimate to look for (and we also gave an elementary proof of this fact in the spirit of E. Gagliardo and L. Nirenberg's slicing argument).    However, as we will see in what follows, this suggests many more open questions to be addressed.  The purpose of this paper is to give an introduction and exposition of the author's perspective of this area and to discuss some open problems in this regard.

The remainder of the paper is dedicated to discussing the connections of the material presented in this section with various literature on the topic, to provide some proofs of the inequalities in the introduction, and to prepare the reader for the open problems in the last section.  In particular, Section \ref{potentials} presents a curious formula for the logarithm which was obtained in collaboration with Rahul Garg in the papers \cite{GargSpector,GargSpector1} and its relation to the the recent work of J. Bourgain and H. Brezis \cite{Bourgain_Brezis_2002, Bourgain_Brezis_2003, Bourgain_Brezis_2004, Bourgain_Brezis_2007}, as well as the more classical work of F. John and L. Nirenberg \cite{JohnNirenberg}.  In Section \ref{potentialestimates} we give some proofs of the results which emphasizes the connections with work of L. Hedberg \cite{Hedberg}, S. Sobolev \cite{sobolev} and A. Zygmund \cite{Zygmund}.  In Section \ref{vector} we discuss the case of vector-valued inequalities, where certain algebraic conditions become relevant in the determination of whether a given differential operator can support a Sobolev inequality.  In particular, while one has a characterization of the differential operators which yield an embedding into the Lebesgue spaces - the elliptic and canceling operators of J. Van Schaftingen \cite{VanSchaftingen_2013}, whether these operators support the improvements known in the classical setting has still not been resolved.  In Section \ref{improvements} we present some results in Lorentz spaces which are the optimal known inequalities for two settings:  estimates for Riesz potentials due to the author and estimates for elliptic and $(d-1)-$canceling operators due to the author and J. Van Schaftingen.  In Section \ref{traceinequalities} we discuss the trace inequality of N. Meyers and W.P. Ziemer and how it represents the best known Sobolev inequality in this classical first order setting.  Finally, in Section \ref{openproblems} we introduce and discuss some open problems the author feels would yield some insight into this question of Sobolev inequalities.

\section{A Curious Formula for the Logarithm (Connection with the work of J. Bourgain and H. Brezis)} \label{potentials}

One has a number of avenues to derive a formula for the solution of \eqref{poisson}, for example by studying an ordinary differential equation or utilizing the Fourier transform.  Whatever the method, if $d\geq 3$ one can verify that
\begin{align}
u(x):= \frac{1}{(d-2)|S^{d-1}|} \int_{\mathbb{R}^d} \frac{f(y)}{|x-y|^{d-2}}\;dy \label{newton}
\end{align}
satisfies \eqref{poisson} in the sense of distributions, i.e.
\begin{align}
-\int_{\mathbb{R}^d} u \Delta \varphi\;dx = \int_{\mathbb{R}^d} f\varphi \;dx \label{dist_poisson}
\end{align}
for all $\varphi \in C^\infty_c(\mathbb{R}^d)$.  In the case $d=2$ one finds the solution to \eqref{poisson} in the appropriate sense, i.e. \eqref{dist_poisson} is given by the logarithmic potential
\begin{align}
u(x):= \frac{1}{2\pi} \int_{\mathbb{R}^2} \log \frac{1}{|x-y|} f(y)\;dy. \label{logarithmic}
\end{align}

However, recently in \cite{GargSpector,GargSpector1} Rahul Garg and the author gave a new representation of the solution \eqref{logarithmic} that does not involve the Logarithm.  Our motivation for doing so stems from the fact that one when $f \in L^p(\mathbb{R}^2)$ and $1<p\leq2$ one expects continuity estimates, with the case $p=2$ corresponding to the almost Lipschitz estimate of H. Brezis and S. Wainger \cite{BrezisWainger}.  Yet the standard approach to continuity estimates for Riesz potentials does not apply in such a setting.   This led us to the following approach.

Let us denote by 
\begin{align*}
\widehat{\varphi}(\xi) = \int_{\mathbb{R}^d} \varphi(x) e^{-2\pi i x\cdot \xi}\;dx
\end{align*}
the Fourier transform of a function $\varphi :\mathbb{R}^d \to \mathbb{R}$.  Then for sufficiently nice $u$ which satisfies \eqref{poisson} in the sense of distributions, we must have
\begin{align}
\widehat{u}(\xi) = \frac{\widehat{f}(\xi)}{(2\pi |\xi|)^2}.\label{fouriertransformsolution}
\end{align}
Thus for $d \geq 3$, as the Fourier transform takes convolution to multiplication and vice versa, we conclude as usual that
\begin{align*}
u(x)=I_2\ast f(x) = I_2f(x),
\end{align*}
which is our reason for suggestively writing \eqref{fouriertransformsolution}, since for general $\alpha \in (0,d)$ one has
\begin{align*}
\widehat{I_\alpha f}(\xi) = \frac{\widehat{f}(\xi)}{(2\pi |\xi|)^\alpha}.
\end{align*}
In fact, this determines precisely the constant
\begin{align*}
\gamma(\alpha):= \frac{\pi^{d/2}2^\alpha \Gamma\left(\frac{\alpha}{2}\right)}{\Gamma\left(\frac{d-\alpha}{2}\right)}.
\end{align*}

In the case $d=2$ the denominator of the preceding equation \eqref{fouriertransformsolution} is not locally integrable near zero.  This is not a problem, as it forces one to impose $\hat{f}(0)=0$ to correct this, which is simply $\int f=0$, and then a suitable limiting process allows one to obtain \eqref{logarithmic}.  However, the appearance of the logarithm was not amenable to the estimates Rahul and the author wanted to show, which led us to the factorization
\begin{align*}
\widehat{u}(\xi) = \frac{-1}{4\pi^2}\frac{i \xi}{|\xi|^3}  \cdot   \frac{i\xi}{|\xi|} \widehat{f}(\xi)
\end{align*}
Now one can check that 
\begin{align*}
\widehat{R f}(\xi) = \frac{i\xi}{|\xi|} \widehat{f}(\xi),
\end{align*}
while in the case $d=2$ the Riesz transform has kernel $\frac{x}{|x|^3}$, up to a multiplicative constant.  In particular the scaling suggests that multiplication by the term
\begin{align*}
\frac{i \xi}{|\xi|^3}
\end{align*}
in Fourier spaces should invert as convolution with the bounded function
\begin{align*}
x \mapsto c\frac{x}{|x|}
\end{align*}
for some appropriate constant $c$.  Indeed, this can be made precise to yield
\begin{align*}
I_2f(x) = \frac{1}{2\pi}\int_{\mathbb{R}^2} \frac{x-y}{|x-y|} \cdot Rf(y)\;dy,
\end{align*}
which gives a new representation of the fundamental solution to Poisson's equation in the plane.  

This idea generalizes to the logarithmic potential in any number of dimensions, and we here recall the computation from the recent paper of the author and Itai Shafrir \cite{Shafrir-Spector} which works out the details.  In particular we instead rely on the semi-group property of the Riesz potentials to write
\begin{align*}
I_d f=  I_{d-1} I_1f &= \frac{1}{\gamma(d-1)} \int_{\mathbb{R}^d} I_1f(y) \frac{1}{|x-y|}\;dy \\
&=\frac{1}{(d-1)\gamma(d-1)} \int_{\mathbb{R}^N} I_1f(y) \operatorname*{div} \left(\frac{x-y}{|x-y|}\right)\;dy.
\end{align*}
However, now performing an integration by parts we have
\begin{align*}
 I_d f(x)= \frac{1}{(d-1)\gamma(d-1)} \int_{\mathbb{R}^d} Rf(y) \cdot  \frac{x-y}{|x-y|}\,dy.
\end{align*}

This formula sheds some light on classical results in the theory of Hardy spaces and $BMO$ that we now discuss.  From the theory of Hardy spaces of E. Stein and G. Weiss, the Hardy space consists of functions $f \in L^1(\mathbb{R}^d)$ such that their Riesz transform $Rf \in L^1(\mathbb{R}^d;\mathbb{R}^d)$.  Thus the duality of $\mathcal{H}^1(\mathbb{R}^d)$ and the John-Nirenberg space of functions of bounded mean oscillation ($BMO$) obtained by C. Fefferman \cite{Fefferman} (see also Fefferman and Stein \cite{FeffermanStein}) implies that every $g \in BMO(\mathbb{R}^d)$ can be expressed as
\begin{align*}
g = g_0 + \sum_{j=1}^d R_jg_j
\end{align*}
for some $\{g_j\}_{j=0}^d \subset L^\infty(\mathbb{R}^d)$ (a constructive proof of this fact was subsequently obtained by A. Uchiyama \cite{Uchiyama}).  

The work of J. Bourgain and H. Brezis \cite{Bourgain_Brezis_2002,Bourgain_Brezis_2003}, among other results, demonstrates an improvement to this representation for certain $BMO$ functions, those of the form $I_1f$ for $f \in L^d(\mathbb{R}^d)$, as they show that such functions have a representation
\begin{align*}
I_1f = \sum_{j=1}^d R_jg_j
\end{align*}
for some $\{g_j\}_{j=1}^d \subset L^\infty(\mathbb{R}^d)$, that is, one does not need the function $g_0$.  

It is at this point that our result \eqref{potentialnodirac} enters, since it implies, by duality, that such a representation extends to any $\alpha \in (0,d)$, i.e. for all $f \in L^{d/\alpha}(\mathbb{R}^d)$ one has
\begin{align*}
I_\alpha f = \sum_{j=1}^d R_jg_j
\end{align*}
for some $\{g_j\}_{j=1}^d \subset L^\infty(\mathbb{R}^d)$.

In general these functions $\{g_j\}_{j=1}^d \subset L^\infty(\mathbb{R}^d)$ are not explicit and cannot be chosen linearly.  However, the preceding calculation taken from the papers \cite{GargSpector,GargSpector1,Shafrir-Spector} show that for the canonical example of a $BMO$ function, $\log |x|$, one does not need $g_0$ and explicitly
\begin{align*}
\log |x| = \sum_{j=1}^d R_j  \frac{1}{(d-1) \gamma(d-1)}\frac{x_j}{|x|}.
\end{align*}
In fact, returning to the seminal paper \cite{JohnNirenberg}, F. John and L. Nirenberg had given a family of examples of functions of bounded mean oscillation in a bounded domain $\Omega \subset \mathbb{R}^d$, 
\begin{align*}
u(x) = \int_\Omega \log |x-y| f(y)\;dy.
\end{align*}
In particular, one follows the above result to see that all of these functions are special functions in $BMO$ which can be expressed as
\begin{align*}
g =\sum_{j=1}^d R_jg_j,
\end{align*}
with explicit $\{g_j\}_{j=1}^d \subset L^\infty(\mathbb{R}^d)$ given by
\begin{align*}
g_j(x) := \frac{1}{(d-1) \gamma(d-1)} \int_{\Omega} \frac{x_j-y_j}{|x-y|} f(y)\;dy.
\end{align*}

\section{A Few Proofs (A Unified Approach to S. Sobolev and A. Zygmund's theorems)}\label{potentialestimates}

It was observed by L. Hedberg in \cite{Hedberg} that one can give a proof of S. Sobolev's theorem for the Riesz potentials (here recorded as Theorem \ref{Sobolev}) by a simple pointwise estimate.  We first show here how this pointwise estimate yields S. Sobolev's theorem and A. Zygmund's weak-type estimate for the Riesz potentials before returning to prove some estimates for Poisson's equation.

Let us therefore recall the pointwise inequality of L. Hedberg \cite{Hedberg}:  If $1\leq p < d/\alpha$, then
\begin{align}
|I_\alpha f(x)| \leq C \mathcal{M}(f)(x)^{1-\alpha p/d} \|f\|_{L^p(\mathbb{R}^d)}^{\alpha p /d}\label{Hedbergest}
\end{align}
for all $f \in L^p(\mathbb{R}^d)$.  Here we utilize the notation $\mathcal{M}(f)$ to denote the Hardy-Littlewood maximal function of a function $f$ by
\begin{align*}
\mathcal{M}(f)(x):= \sup_{r>0} \frac{1}{|B(x,r)|} \int_{B(x,r)} |f(y)|\;dy,
\end{align*}
for which we require the inequalities
\begin{align*}
\sup_{t>0} t|\{ |\mathcal{M}f(x)>t\}| &\leq C \|f\|_{L^1(\mathbb{R}^d)} \\
\|\mathcal{M}f\|_{L^p(\mathbb{R}^d)} &\leq C' \|f\|_{L^p(\mathbb{R}^d)} \quad 1<p \leq +\infty.
\end{align*}
These estimates are standard, see e.g. \cite{Sharmonic}.

Note that from this one obtains immediately
\begin{align*}
\|I_\alpha f\|_{L^q(\mathbb{R}^d)} &\leq C \|  \mathcal{M}(f)^{1-\alpha p/d} \|f\|_{L^p(\mathbb{R}^d)}^{\alpha p /d}\|_{L^q(\mathbb{R}^d)} \\
&= C\|  \mathcal{M}(f)\|_{L^p(\mathbb{R}^d)}^{1-\alpha p/d} \|f\|_{L^p(\mathbb{R}^d)}^{\alpha p /d} \\
&\leq C'\|f\|_{L^p(\mathbb{R}^d)}
\end{align*}
by the boundedness of the maximal function and combining like terms, which is Theorem \ref{Sobolev}.  However, the same proof using the weak-$(1,1)$ estimate for the maximal function implies
\begin{theorem}
Let $\alpha \in (0,d)$.  Then there exists a constant $C=C(\alpha,d)>0$ such that
\begin{align*}
\sup_{t>0} t |\{ |I_\alpha f (x)|>t \}|^\frac{d-\alpha}{d} \leq C\|f\|_{L^1(\mathbb{R}^d)}
\end{align*}
for all $f \in L^1(\mathbb{R}^d)$.
\end{theorem}
We provide the details for the convenience of the reader.
\begin{proof}
\begin{align*}
\{ |I_\alpha f (x)|>t \} \subset \{ C\mathcal{M}f(x)^{1-\alpha/d} \|f\|^{\alpha/d}_{L^1(\mathbb{R}^d)}>t \},
\end{align*}
which yields for every $t>0$ the estimate
\begin{align*}
t |\{ |I_\alpha f (x)|>t \}|^\frac{d-\alpha}{d} \leq t| \{ C\mathcal{M}f(x)^{1-\alpha/d} \|f\|_{L^1(\mathbb{R}^d)}^{\alpha/d}>t \}|^\frac{d-\alpha}{d}.
\end{align*}
We let 
\begin{align*}
s = \left( \frac{t}{C\|f\|_{L^1(\mathbb{R}^d)}^{\alpha/d}}\right)^{d/(d-\alpha)}
\end{align*}
and find
\begin{align*}
t |\{ |I_\alpha f (x)|>t \}|^\frac{d-\alpha}{d}  &\leq t| \{ C\mathcal{M}f(x)^{1-\alpha/d} \|f\|_{L^1(\mathbb{R}^d)}^{\alpha/d}>t \}|^\frac{d-\alpha}{d} \\
&= s^{(d-\alpha)/d} C \|f\|^{\alpha/d}_{L^1(\mathbb{R}^d)} | \{ \mathcal{M}f(x)>s \}|^\frac{d-\alpha}{d} \\
&\leq C \|f\|^{\alpha/d}_{L^1(\mathbb{R}^d)} \left(C' \|f\|_{L^1(\mathbb{R}^d)}\right)^{(d-\alpha)/d}\\
&= C'' \|f\|_{L^1(\mathbb{R}^d)},
\end{align*}
by utilizing the weak-type estimate for $\mathcal{M}$ and combining like terms.  The desired conclusion follows by taking the supremum in $t>0$.
\end{proof}

In particular, if $d\geq 3$, one immediately deduces the inequalities \eqref{u} and \eqref{u'}, while \eqref{gradu} and \eqref{gradu'} can be argued as follows.  We have
\begin{align*}
\nabla u(x) = \nabla I_1 I_1f(x) \equiv R(I_1f)(x),
\end{align*}
and so
\begin{align*}
\|\nabla u \|_{L^q(\mathbb{R}^d)} &= \|R(I_1f) \|_{L^q(\mathbb{R}^d)} \\
&\leq C\|I_1f \|_{L^q(\mathbb{R}^d)} \\
&\leq C' \|f\|_{L^p(\mathbb{R}^d)}
\end{align*}
by the boundedness of the Riesz transforms
\begin{align*}
R_i :L^p(\mathbb{R}^d) \to L^p(\mathbb{R}^d) \text{ for } 1<p<+\infty
\end{align*}
and Theorem \ref{Sobolev} (see \cite{Sharmonic} or \cite{grafakos}).  The same argument applies to obtain \eqref{gradu'}, since the fact that the Riesz transforms are bounded on $L^p(\mathbb{R}^d)$ for $1<p<+\infty$ implies they are bounded on $L^{p,\infty}(\mathbb{R}^d)$ for the same values of $p$ (see \cite{grafakos}).

The last remaining items to discuss in this section are the inequalities \eqref{hessianu} and \eqref{hessianu'}.  The former is straightforward, as
\begin{align*}
-\Delta u \equiv -\frac{\partial^2u}{\partial x_i \partial x_j} = -\frac{\partial^2}{\partial x_i \partial x_j}I_2 f =  -R_i R_j f
\end{align*}
for $R_i,R_j$ the $i$th and $j$th components of the vector valued Riesz transform $R$, respectively and so the result follows again by boundedness of $R_i,R_j$ on $L^p(\mathbb{R}^d)$.  However, for the weak-type estimate one requires the fact that not only are $R_i$ and $R_j$ Calder\'on-Zygmund operators, but even their composition 
\begin{align*}
Tf:=R_iR_jf
\end{align*}
is such an operator.  Then Lemma 2 of \cite{CalderonZygmund} implies the desired weak-type bound. 

Finally it remains to handle the case $d=2$.  Here the estimates \eqref{gradu} and \eqref{hessianu} can be argued in a similar manner, and so we are left to discuss  the estimates \eqref{u} and \eqref{u'}.  But here one has only the case $p=1$ and so we are reduced to proving some type of replacement for \eqref{u'} (since $d-2=0$).  Yet the analysis in the preceding section shows
\begin{align*}
u(x) = \frac{1}{\gamma(1)}\sum_{i=1}^2 R_i \int_{\mathbb{R}^2} \frac{x_i-y_i}{|x-y|} f(y)\;dy,
\end{align*}
which shows $u \in BMO(\mathbb{R}^d)$.

It is interesting to note that, as in the introduction, we replace $\|f\|_{L^1(\mathbb{R}^2)}$ with $\|Rf\|_{L^1(\mathbb{R}^2;\mathbb{R}^2)}$ we obtain the inequality
\begin{align*}
\|I_2 f\|_{L^\infty(\mathbb{R}^2)} \leq C \|Rf\|_{L^1(\mathbb{R}^2;\mathbb{R}^2)},
\end{align*}
and what is more, we also have
\begin{align*}
\| I_2 Rf\|_{L^\infty(\mathbb{R}^d)} \leq C \|I_1 f\|_{L^{2,1}(\mathbb{R}^2)}\leq C' \|Rf\|_{L^1(\mathbb{R}^2;\mathbb{R}^2)}.
\end{align*}
Here we have utilized a Lorentz space extension of potential embeddings, see Section \ref{improvements} for a further excursion into this area.  That is, if $Rf \in L^1(\mathbb{R}^2;\mathbb{R}^2)$, then both $I_2 f$ and its Riesz transform $R I_2 f$ are bounded functions.  The same is true for $I_df$, $RI_df$ if we assume $Rf \in L^1(\mathbb{R}^d;\mathbb{R}^d)$.

\section{Vector Inequalities (Connection with the work of J. Van Schaftingen)}\label{vector}

A vector analogue of the classical inequality of E. Gagliardo \cite{Gagliardo} and L. Nirenberg \cite{Nirenberg} follows easily from the same argument, yet such an inequality is not optimal, as one does not need the full gradient in order to obtain an embedding into $L^{d/(d-1)}(\mathbb{R}^d)$.  For example, one has the 1973 result of M.J. Strauss \cite{Strauss}, which asserts that one only needs the symmetric part of the gradient when $u :\mathbb{R}^d \to \mathbb{R}^d$.  In particular if one defines the linearized deformation tensor
\begin{align*}
    Eu:= \frac{1}{2}\bigl(\nabla u + (\nabla u)^T\bigr),
\end{align*}
then M.J. Strauss proves the inequality
\begin{align}\label{Strauss}
\|u \|_{L^{d/(d-1)}(\mathbb{R}^d;\mathbb{R}^d)} \leq C \| Eu \|_{L^1(\mathbb{R}^d,\mathbb{R}^{d\times d})}
\end{align}
for all sufficiently nice $u$.  More generally, J. Van Schaftingen has shown that the vector differential inequality
\begin{align}
\label{ineq_EC}
\|u \|_{L^{d/(d-1)}(\mathbb{R}^d;V)} \leq C \| A (D) u \|_{L^1(\mathbb{R}^d,E)}
\end{align}
holds for every vector field \(u \in C^\infty_c (\mathbb{R}^d, V)\) 
if and only if the homogeneous first-order linear vector differential operator with constant coefficients \(A (D)\) is elliptic and canceling \cite{VanSchaftingen_2013}*{Theorem 1.3}.  Here we recall the definition of a canceling operator.

\begin{definition}
Let \(\ell \in \{0, \dotsc, d\}\).
A homogeneous differential operator with constant coefficients \(A (D)\) is \emph{ \(\ell\)--canceling} whenever 
\begin{align*}
    \bigcap_{\substack{W \subseteq \mathbb{R}^d\\ \dim W = \ell}}
      \linspan 
        \,
        \bigl\{ 
          A (\xi)[v] 
        \;:\; 
          \xi \in W 
          \text{ and } 
          v \in V
        \bigr\}
  =
    \{0\}
  .
\end{align*}
\end{definition}

The term canceling is precisely \(1\)--canceling as defined in \cite{VanSchaftingen_2013}*{Definition 1.2}.

One then wonders whether in the vector-valued setting the Riesz potentials admit similar improvements.  Indeed this is the case, as from the work of Van Schaftingen in \cite{VanSchaftingen_2013}, with the argument of \cite{SSVS}, one obtains 
\begin{theorem}
Let $\alpha \in (0,d)$, $V, E$ be finite dimensional Banach spaces, and suppose  the homogeneous first-order linear differential operator with constant coefficients $A:C^\infty_c(\mathbb{R}^d;V) \to C^\infty_c(\mathbb{R}^d;E)$ is elliptic and canceling.  Then there exists a constant $C>0$ such that
\begin{align*}
\|I_\alpha f \|_{L^{d/(d-\alpha)}(\mathbb{R}^d;V)} \leq C \|A(D)I_1f\|_{L^1(\mathbb{R}^d;E)}
\end{align*}
for all $f \in C^\infty_c(\mathbb{R}^d;V)$.
\end{theorem}
In particular one can apply Theorem 8.3 in \cite{VanSchaftingen_2013} to obtain such a result by recognizing that the Triebel-Lizorkin spaces coincide with the above space of Riesz potentials of $L^{d/(d-\alpha)}$ functions.   A similar statement would also seem to hold in the case of higher order homogeneous linear differential operator with constant coefficients.  Thus, this gives the complete picture for the mapping properties of Riesz potentials into Lebesgue spaces in the vector-valued setting.

\section{Some Optimal Lorentz Space Estimates}\label{improvements}
In this section we discuss improvements to the preceding embeddings on the Lorentz scale.  The Lorentz spaces arise naturally in the study of PDE and are of interest for a number of reasons:  They appear immediately in approximation theory and the interpolation of Banach spaces when one only begins with an interest in the classical Lebesgue spaces; they are a scale of spaces which incorporate the natural weak-type spaces in harmonic analysis; the extension of H\"older and Young's inequality to this scale imply Theorem \ref{Sobolev} (see Theorem \ref{Lorentz} below); they can be used to establish existence of solutions to the wave and Schr\"odinger equation (see, e.g. Keel and Tao \cite{KeelTao}), and can be used to deduce continuity and Lipschitz continuity in the context of Harmonic maps to manifolds (see H\'elein \cite{Helein}).  The value of such spaces has been known to experts for some time - the author was aware of many results from L. Tartar's article \cite{tartar} from 1998, while in the lecture of this material at Rutgers, H. Brezis referred him to his article on the subject \cite{BrezisLorentz} from 1982.  Let us also mention two further references which may be of interest to the reader, the classical papers of R. Hunt \cite{Hunt} and R. O'Neil \cite{oneil}.  

In the literature there are a number of possible definitions of the Lorentz spaces $L^{p,q}(\mathbb{R}^d)$.  We find it most useful for the unacquainted reader to see these spaces as a refinement of the Lebesgue spaces which can be simply understood in the following manner.   By Cavalieri's principle, one begins by expressing the $L^p$-norm as integration of the superlevel sets
\begin{align}
\|f\|_{L^p(\mathbb{R}^d)}^p = p \int_0^\infty t^p | \{ |f|>t\}|\;\frac{dt}{t}. \label{cavaliere}
\end{align}
This shows that belonging to $L^p$ requires specific decay such that the function 
\[t \mapsto (t | \{ |f|>t\}|^\frac{1}{p})^{p}\]
is integrable near both $t=0$ and $t=+\infty$ with respect to the measure $\frac{dt}{t}$, which requires that it tends to zero at both endpoints.  Therefore, any power of this map must tend to zero, though is not necessarily integrable.  Varying the exponent here as a second parameter $r$ and imposing integrability, one obtains the Lorentz spaces $L^{p,r}(\mathbb{R}^d)$, with quasi-norm for $1<p<+\infty$ and $1\leq r<+\infty$
\begin{align}
\|f\|_{L^{p,r}(\mathbb{R}^d)}^r := p \int_0^\infty (t | \{ |f|>t\}|^\frac{1}{p})^r\;\frac{dt}{t} \label{Lorentzquasinorm}
\end{align}
and $1\leq p <+\infty$ and $r=+\infty$
\begin{align*}
\|f\|_{L^{p,\infty}(\mathbb{R}^d)} := \sup_{t>0}  t | \{ |f|>t\}|^\frac{1}{p}.
\end{align*}
One then observes that $L^{p,p}(\mathbb{R}^d) \equiv L^p(\mathbb{R}^d)$, while more generally one obtains a scale of spaces which is nested with respect to the second parameter, i.e.
\begin{align*}
L^{p,1}(\mathbb{R}^d) \subset L^{p,p}(\mathbb{R}^d) = L^p(\mathbb{R}^d) \subset  L^{p,\infty}(\mathbb{R}^d).
\end{align*}
For these spaces, one has the analogue of H\"older's inequality (Theorem 3.4 in \cite{oneil}):
\begin{theorem}\label{holder}
Let $f \in L^{q_1,r_1}(\mathbb{R}^d)$ and $g \in L^{q_2,r_2}(\mathbb{R}^d)$, where
\begin{align*}
\frac{1}{q_1}+\frac{1}{q_2}&=\frac{1}{q}<1\\
\frac{1}{r_1}+\frac{1}{r_2}&\geq  \frac{1}{r},
\end{align*}
for some $r \geq 1$.   Then
\begin{align*}
\|fg\|_{L^{q,r}(\mathbb{R}^d)} \leq q'\|f \|_{L^{q_1,r_1}(\mathbb{R}^d)}\|g \|_{L^{q_2,r_2}(\mathbb{R}^d)}
\end{align*}
\end{theorem}

One also has Young's inequality (Theorem 3.1 in \cite{oneil}):
\begin{theorem}\label{young}
Let $f \in L^{q_1,r_1}(\mathbb{R}^d)$ and $g \in L^{q_2,r_2}(\mathbb{R}^d)$, and suppose $1< q<+\infty$ and $1\leq r\leq +\infty$ satisfy
\begin{align*}
\frac{1}{q_1}+\frac{1}{q_2}-1&=\frac{1}{q}\\
\frac{1}{r_1}+\frac{1}{r_2}&\geq \frac{1}{r}.
\end{align*}
Then
\begin{align*}
\|f\ast g\|_{L^{q,r}(\mathbb{R}^d)} \leq 3q \|f \|_{L^{q_1,r_1}(\mathbb{R}^d)}\|g \|_{L^{q_2,r_2}(\mathbb{R}^d)}.
\end{align*}
\end{theorem}
Now $I_\alpha \notin L^r(\mathbb{R}^d)$ for $1\leq r \leq +\infty$, while it is an exercise to show that $I_\alpha \in L^{d/(d-\alpha),\infty}(\mathbb{R}^d)$.  Thus this inequality implies the following improvement to Theorem \ref{Sobolev}.
\begin{theorem} \label{Lorentz}
Let $0<\alpha<d$ and $1<p<d/\alpha$.  Then there exists a constant 
$C=C(p,\alpha,d)>0$ such that
\begin{align}\label{oneil}
\|I_\alpha f \|_{L^{q,p}(\mathbb{R}^d)} \leq C \|f\|_{L^p(\mathbb{R}^d)}
\end{align}
for all $f \in L^p (\mathbb{R}^d)$, where
\begin{align*}
\frac{1}{q}=\frac{1}{p}-\frac{\alpha}{d}.
\end{align*}  
\end{theorem}

The model inequality in this setting of Lorentz spaces is the result of A. Alvino \cite{Alvino} who proved the inequality
\begin{align}
\label{Alvino}
\|u \|_{L^{d/(d-1),1}(\mathbb{R}^d)} \leq C' \|\nabla u \|_{L^1(\mathbb{R}^d,\mathbb{R}^d)}
\end{align}
holds for all functions $u \in W^{1,1}(\mathbb{R}^d)$.  Here we observe that this improves the inequality of E. Gagliardo \cite{Gagliardo} and L. Nirenberg \cite{Nirenberg}, and so it is natural to ask whether similar improvements can be made for the Riesz potentials.

Indeed, recently in \cite{Spector} the author proved the optimal Lorentz inequality for the Riesz potentials, the following
\begin{theorem}\label{mainresult}
Let $d\geq 2$ and $\alpha \in (0,d)$.  There exists a constant $C=C(\alpha,d)>0$ such that
\begin{align}\label{l1typeestimate}
\| I_\alpha f \|_{L^{d/(d-\alpha),1}(\mathbb{R}^d)} \leq C\| R f \|_{L^1(\mathbb{R}^d;\mathbb{R}^d)}
\end{align}
for all $f \in C^\infty_c(\mathbb{R}^d)$ such that $Rf \in L^1(\mathbb{R}^d;\mathbb{R}^d)$.
\end{theorem}
The key idea of the proof was that one has an equivalent formulation of this inequality as the inequality
\begin{align}\label{gradientinequality}
\| I_\alpha \nabla u \|_{L^{d/(d-\alpha),1}(\mathbb{R}^d;\mathbb{R}^d)} \leq C\| \nabla u \|_{L^1(\mathbb{R}^d;\mathbb{R}^d)},
\end{align}
which can argued by using the fact that $Rf=\nabla I_1 u$ is a gradient along with the boundedness of the Riesz transforms on $L^{p,r}(\mathbb{R}^d)$ for all $1<p<+\infty$, $1\leq r \leq +\infty$.  In such a formulation, the idea of V. Maz'ya \cite{mazya} reduces the question to that of proving an isoperimetric inequality, and this is the main work of the paper \cite{Spector}.  

This prompts one to wonder whether similar improvements can be made for first-order elliptic and canceling operators.  The complete answer to the former question is only known in the plane where \((d - 1)\)--canceling is precisely canceling, as the author and J. Van Schaftingen have shown in \cite{Spector-VanSchaftingen-2018} the following
\begin{theorem}
Let $V$ and $E$ be finite-dimensional spaces and suppose that the homogeneous first-order linear differential operator with constant coefficients \(A (D):C^\infty_c (\mathbb{R}^d, V) \to C^\infty_c (\mathbb{R}^d, E)\) is elliptic and \((d - 1)\)--canceling.  Then there exists a constant $C>0$ such that 
\begin{align*}
\|u \|_{L^{d/(d-1),1}(\mathbb{R}^d;V)} \leq C \| A (D) u \|_{L^1(\mathbb{R}^d,E)}
\end{align*}
for every \(u \in C^\infty_c (\mathbb{R}^d, V)\).
\end{theorem}

Note that one example of a \((d - 1)\)--canceling operator in any number of dimensions is the deformation operator, and so this result contains the optimal Lorentz inequality
\begin{align*}
\|u \|_{L^{d/(d-1),1}(\mathbb{R}^d;\mathbb{R}^d)} \leq C \| Eu \|_{L^1(\mathbb{R}^d,\mathbb{R}^{d\times d})}.
\end{align*}

This leads to further open problems which we discuss in Section \ref{openproblems}.

\section{Trace Inequalities (Connection with the work of N. Meyers and W.P. Ziemer)}\label{traceinequalities}

A 1977 result of Meyers and Ziemer \cite{Meyers-Ziemer-1977} asserts the existence of a constant $C>0$ such that one has the inequality
\begin{align}
\int_{\mathbb{R}^d} |u| \;d\mu \leq C \int_{\mathbb{R}^d} |\nabla u|\;dx \label{meyersziemer}
\end{align}
for every $u \in W^{1,1}(\mathbb{R}^d)$ and every non-negative Radon measure $\mu$ satisfying the ball growth condition
\begin{align*}
\mu(B(x,r)) \leq C'r^{d-1}
\end{align*}
for all $x \in \mathbb{R}^d$, $r>0$, and some $C'>0$.

This inequality represents the state of the art when compared with all of the inequalities discussed so far, which we now show before giving a sketch of its proof.  In particular, let us show how this trace inequality implies\footnote{These connections have been discussed by A. Ponce and the author in \cite{Ponce-Spector}.} the Sobolev inequality of E. Gagliardo \cite{Gagliardo} and L. Nirenberg \cite{Nirenberg}, its Lorentz improvement \cite{Alvino}, and even Hardy's inequality.    To this end, let $g \in L^{d,\infty}(\mathbb{R}^d)$ be non-negative and define, for $A\subset \mathbb{R}^d$ measurable, the non-negative measure
\begin{align*}
\mu(A):=\int_A g(y)\;dy.
\end{align*}
Then H\"older's inequality in the Lorentz spaces implies, for every $B(x,r)\subset \mathbb{R}^d$,
\begin{align*}
\mu(B(x,r)) &= \int_{B(x,r)} g(y)\;dy \leq \| g\|_{L^{d,\infty}(\mathbb{R}^d)} \| \chi_{B(x,r)}\|_{L^{d/(d-1),1}(\mathbb{R}^d)} \\
&= \| g\|_{L^{d,\infty}(\mathbb{R}^d)} |B(x,r)|^{1-1/d}\\
&= \| g\|_{L^{d,\infty}(\mathbb{R}^d)} |B(0,1)|^{1-1/d} r^{d-1}.
\end{align*}
Thus, the trace inequality implies
\begin{align*}
\int_{\mathbb{R}^d} |u(x)| g(x)\;dx \leq C \| g\|_{L^{d,\infty}(\mathbb{R}^d)}  \int_{\mathbb{R}^d} |\nabla u(x)|\;dx,
\end{align*}
and as the norm in $L^{d/(d-1),1}(\mathbb{R}^d)$ can be realized via duality, i.e.
\begin{align*}
\|u\|_{L^{d/(d-1),1}(\mathbb{R}^d)} = \sup_{g \in L^{d,\infty}(\mathbb{R}^d)} \int_{\mathbb{R}^d} |u(x)| \frac{g(x)}{\| g\|_{L^{d,\infty}(\mathbb{R}^d)}}\;dx,
\end{align*}
the Lorentz estimate of A. Alvino \cite{Alvino} follows.  This same argument can be applied in the case
\begin{align*}
\mu(A):=\int_A \frac{1}{|y|}\;dy.
\end{align*}
to obtain Hardy's inequality as a corollary, and
\begin{align*}
\mu(A):=\int_A g(y)\;dy
\end{align*}
for $g \in L^d(\mathbb{R}^d)$ to obtain the inequality of E. Gagliardo \cite{Gagliardo} and L. Nirenberg \cite{Nirenberg}.  

For the convenience of the reader we recall the proof of the inequality \eqref{meyersziemer}.  First, one begins with the Poincar\'e inequality
\begin{align*}
\fint_{B(x,r)} |u(y) - \fint_{B(x,r)} u| \;dy \leq C r \fint_{B(x,r)} |\nabla u(y)|\;dy,
\end{align*}
which can be obtained by the fundamental theorem of calculus (see, e.g. p.~142 in \cite{EvansGariepy}, which extends to functions of bounded variation by density).  Then one takes $u=\chi_E$ for a set $E \subset \mathbb{R}^d$ such that $D\chi_E$ is a Radon measure.  The preceding inequality is written for smooth functions, but extends to $BV_{loc}(\mathbb{R}^d)$.  Thus
\begin{align*}
\fint_{B(x,r)} \fint_{B(x,r)} |\chi_E(y) - \chi_E(z)| \;dy dz \leq 2C r^{1-d} |D\chi_E|(B(x,r)).
\end{align*}
But this says
\begin{align*}
 \frac{|E \cap B(x,r)| \times |E^c \cap B(x,r)|}{|B(x,r)| \times |B(x,r)|} &= \frac{1}{2} \fint_{B(x,r)} \fint_{B(x,r)} |\chi_E(y)-\chi_E(z)|\;dydz \\
& \leq C r^{1-d} |D\chi_E|(B(x,r)).
\end{align*}
Now, for $E \subset \mathbb{R}^d$ open, bounded, and of finite perimeter and $x \in E$, the map 
\begin{align*}
x \mapsto \frac{|E \cap B(x,r)|}{|B(x,r)|}
\end{align*}
is continuous, equal to one for small $r$ and tends to zero as $r \to \infty$.  Thus the intermediate value theorem guarantees an $r=r_x$ such that
\begin{align*}
\frac{|E \cap B(x,r_x)|}{|B(x,r_x)|} = \frac{1}{2}.
\end{align*}
However, for this same value $r_x$, by finite additivity of the Lebesgue measure, we also have
\begin{align*}
\frac{|E^c \cap B(x,r_x)|}{|B(x,r_x)|} = \frac{1}{2}.
\end{align*}
Therefore we have found an $r_x$ such that
\begin{align}
r_x^{d-1} \leq 4C  |D\chi_E|(B(x,r)). \label{radius}
\end{align}
We are now prepared to estimate
\begin{align*}
\int_{\mathbb{R}^d} |u| \;d\mu,
\end{align*}
and we suppose here that $u \in C^\infty_c(\mathbb{R}^d)$.  The result for more general $u$ can be argued by using maximal function bounds with respect to these Choquet integrals, cf. \cite{Adams:1988}.  In particular we recall the definition of
\begin{align*}
\int_{\mathbb{R}^d} |u| \;d\mu = \int_0^\infty \mu (E_t)\;dt,
\end{align*}
where $E_t:= \{x : |u|>t\}$.  

We now apply the preceding argument to the set $E_t$ to find for each $x \in E_t$ an $r_x$ such that \eqref{radius} holds (and note that $x$ and in turn $r_x$ depend implicitly on $t$).  By Vitali's covering theorem (see, for example, Theorem 1 on p.~27 of \cite{EvansGariepy}) we can find a countable family $\{x_i\} \subset E_t$ and $\{r_i\}$ such that
\begin{align*}
E_t \subset \bigcup_{i=1}^\infty \overline{B(x_i,5r_i)},
\end{align*}
$\overline{B(x_i,r_i)} \cap \overline{B(x_j,r_j)}=\emptyset$ if $i\neq j$ and \eqref{radius} holds.  Therefore
\begin{align*}
\mu(E_t) &\leq \mu \left(\bigcup_{i=1}^\infty \overline{B(x_i,5r_i)}\right) \\
&\leq \sum_{i=1}^\infty C(5r_i)^{d-1} \\
&\leq 5^{d-1} C \sum_{i=1}^\infty |D\chi_{E_t}|(B(x_i,r_i))\\
&\leq 5^{d-1} C |D\chi_{E_t}|(\mathbb{R}^d).
\end{align*}
However, now combining this with the preceding equality we find
\begin{align*}
\int_{\mathbb{R}^d} |u| \;d\mu &= \int_0^\infty \mu(E_t)\;dt \\
&\leq 5^{d-1} C\int_0^\infty |D\chi_{E_t}|(\mathbb{R}^d)\;dt \\
&=5^{d-1} C \int_{\mathbb{R}^d} | \nabla u|\;dx,
\end{align*}
where the last equality follows from the coarea formula.

\section{Open Problems}\label{openproblems}
We are now ready to discuss some open problems in this area, and let us begin with the scalar setting.  In the introduction we saw how the inequality of E. Gagliardo \cite{Gagliardo} and L. Nirenberg \cite{Nirenberg}
\begin{align*}
\|u\|_{L^{d/(d-1)}(\mathbb{R}^d)} \leq C \|\nabla u \|_{L^1(\mathbb{R}^d;\mathbb{R}^d)},
\end{align*}
leads one to predict the potential estimate observed by A. Schikorra, the author, and J. Van Schaftingen: 
\begin{align*}
\|I_\alpha f\|_{L^{d/(d-\alpha)}(\mathbb{R}^d)} \leq C \|R f \|_{L^1(\mathbb{R}^d;\mathbb{R}^d)}.
\end{align*}
Similarly in Section \ref{improvements} we saw how the refinement of A. Alvino \cite{Alvino} on the Lorentz scale,
\begin{align*}
\|u\|_{L^{d/(d-1),1}(\mathbb{R}^d)} \leq C \|\nabla u \|_{L^1(\mathbb{R}^d;\mathbb{R}^d)},
\end{align*}
leads one to predict refinement obtained by the author in \cite{Spector}:
\begin{align*}
\|I_\alpha f\|_{L^{d/(d-\alpha),1}(\mathbb{R}^d)} \leq C(\alpha) \|R f \|_{L^1(\mathbb{R}^d;\mathbb{R}^d)}.
\end{align*}

Yet in the first order setting, we saw that the strongest possible inequality is the trace inequality
\begin{align*}
\int_{\mathbb{R}^d} |u| \;d\mu \leq C \int_{\mathbb{R}^d} |\nabla u(x)|\;dx
\end{align*} 
for all sufficiently nice $u$ and all non-negative Radon measures $\mu$ satisfying the ball growth condition $\mu(B(x,r)) \leq C'r^{d-1}$ for all $B(x,r) \subset \mathbb{R}^d$.

This would seem to motivate
\begin{question}
Let $\alpha \in (0,1)$.  Does one have the existence of a constant $C=C(\alpha,d)>0$ such that
\begin{align*}
\int_{\mathbb{R}^d} |I_\alpha f| \;d\mu \leq C \int_{\mathbb{R}^d} |Rf(x)|\;dx
\end{align*} 
for all sufficiently nice $f$ and all non-negative Radon measures $\mu$ satisfying the ball growth condition $\mu(B(x,r)) \leq C'r^{d-\alpha}$ for all $B(x,r) \subset \mathbb{R}^d$?
\end{question}
However, one has a fundamental obstruction to such an inequality, which has been discussed in a recent work of the author \cite{Spector-PM}.  In particular we there show the impossibility of such an inequality, establishing a crucial difference between the case $\alpha=1$ and $\alpha \in (0,1)$.

While there is no possibility to obtain a trace inequality for the Riesz potentials for $\alpha \in (0,1)$, one wonders
\begin{openproblem}
Can one show the inequality
\begin{align*}
\int_{\mathbb{R}^d} |u| \;d\mu \leq C \int_{\mathbb{R}^d} |\nabla u(x)|\;dx
\end{align*} 
for all sufficiently nice $u$ and all non-negative Radon measures $\mu$ satisfying the ball growth condition $\mu(B(x,r)) \leq C'r^{d-1}$ for all $B(x,r) \subset \mathbb{R}^d$ without using the coarea formula?
\end{openproblem}

In particular such a result may also be of use in the vector setting, which we now discuss.  The work of J. Van Schaftingen shows that the inequality
\begin{align*}
\|u \|_{L^{d/(d-1)}(\mathbb{R}^d;V)} \leq C \|A(D)u\|_{L^1(\mathbb{R}^d;E)}
\end{align*}
holds whenever $A(D)$ is elliptic and canceling.  A Lorentz improvement for this inequality been shown by the author and J. Van Schaftingen in \cite{Spector-VanSchaftingen-2018}:
\begin{align*}
\|u \|_{L^{d/(d-1),1}(\mathbb{R}^d;V)} \leq C \|A(D)u\|_{L^1(\mathbb{R}^d;E)},
\end{align*}
whenever $A(D)$ is elliptic and $(d-1)-$canceling.  This gives a complete resolution to the question in two dimensions, though prompts one to ask
\begin{openproblem}
Let \(d \ge 3\) and let $V$ and $E$ be finite-dimensional spaces.  Further suppose that the first-order homogeneous linear differential operator  with constant coefficients \(A (D):C^\infty_c (\mathbb{R}^d, V) \to C^\infty_c (\mathbb{R}^d, E)\) is elliptic and canceling.  Can one show the existence of a constant $C>0$ such that the inequality
\begin{align*}
\|u \|_{L^{d/(d-1),1}(\mathbb{R}^d;V)} \leq C \|A(D)u\|_{L^1(\mathbb{R}^d;E)},
\end{align*}
holds for every \(u \in C^\infty_c (\mathbb{R}^d, V)\)?
\end{openproblem}
This is Question 1.1 in \cite{Spector-VanSchaftingen-2018}, see also \citelist{\cite{VanSchaftingen_2013}*{Open problem 8.3}\cite{VanSchaftingen_2015}*{Open problem 2}\cite{Bourgain_Brezis_2007}*{Open problem 1}
\cite{VanSchaftingen_2010}*{Open problem 2}}.  A positive answer to this question would give a model inequality for potential estimates in the vector setting and would then lead one to attack the more difficult
\begin{openproblem}
Let \(d \ge 2\) and suppose that the first-order homogeneous linear differential operator  with constant coefficients  \(A (D):C^\infty_c (\mathbb{R}^d, V) \to C^\infty_c (\mathbb{R}^d, E)\) is elliptic and canceling.  Can one show the existence of a constant $C>0$ such that the inequality
\begin{align*}
\|I_\alpha f \|_{L^{d/(d-\alpha),1}(\mathbb{R}^d;V)} \leq C \|A(D)I_1 f\|_{L^1(\mathbb{R}^d;E)},
\end{align*}
holds for every \(f \in C^\infty_c (\mathbb{R}^d, V)\)?
\end{openproblem}

If one is able to establish Lorentz inequalities for vector differential operators, independent of whether the analogous inequality can be established for potentials, a natural question is whether one has the analogue of N. Meyers and W.P. Ziemer's inequality for such operators, perhaps under more restrictive conditions.  The need to make further assumptions on the operator $A(D)$ has been shown in the forthcoming work of F. Gmeineder, B. Rai\cb{t}\v{a}, and J. Van Schaftingen (see  Theorem 1.3 in \cite{GRV}), where they show that the trace inequality
\begin{align*}
\int_{\Sigma} |u| \;d\mathcal{H}^{d-1} \leq C \|A(D)u\|_{L^1(\mathbb{R}^d;E)}
\end{align*}
on any hyperplane $\Sigma \subset \mathbb{R}^d$ is equivalent to the $\mathbb{C}$-ellipticity of $A(D)$ (see also on the earlier work of D. Breit, L. Diening, and F. Gmeineder in \cite{BDG} on the equivalence of $\mathbb{C}$-ellipticity and boundary trace inequalities).  Here we recall that for a first-order homogeneous linear differential operator with constant coefficients \(A (D):C^\infty_c (\mathbb{R}^d, \mathbb{R}^N) \to C^\infty_c (\mathbb{R}^d, \mathbb{R}^K)\), $A(D)$ is  elliptic ($\mathbb{R}$-elliptic) if
\begin{align*}
A(\xi) : \mathbb{R}^N \mapsto \mathbb{R}^K
\end{align*}
is injective for all $\xi \in \mathbb{R}^d\setminus \{0\}$, while $A(D)$ is  $\mathbb{C}$-elliptic\footnote{It has been pointed out to us by Bogdan Rai\cb{t}\v{a} that the terminology of $\mathbb{C}$-ellipticity is due to D. Breit, L. Diening, and F. Gmeineder in \cite{BDG}, while the condition appears at least as early as the references of Smith \cite{Smith,Smith1} and is inspired by Aronszajn \cite{Aronszajn}.} if
\begin{align*}
A(\xi) : \mathbb{C}^N \mapsto \mathbb{C}^K
\end{align*}
is injective for all $\xi \in \mathbb{C}^d\setminus \{0\}$.  Let us remark that $\mathbb{C}$-ellipticity implies canceling (see Lemma 3.2 in \cite{GR}).

This leads us to ask
\begin{openproblem}\label{traceEC}
Let \(d \ge 2\) and suppose that the first-order homogeneous linear differential operator  with constant coefficients \(A (D):C^\infty_c (\mathbb{R}^d, V) \to C^\infty_c (\mathbb{R}^d, E)\) is $\mathbb{C}$-elliptic.  Can one show the existence of a constant $C>0$ such that the inequality
\begin{align*}
\int_{\mathbb{R}^d} |u| \;d\mu  \leq C \|A(D)u\|_{L^1(\mathbb{R}^d;E)},
\end{align*}
holds for every \(u \in C^\infty_c (\mathbb{R}^d, V)\) and every non-negative Radon measure $\mu$ such that $\mu(B(x,r)) \leq C'r^{d-1}$ for all $B(x,r) \subset \mathbb{R}^d$?
\end{openproblem}

A partial answer in the case of hyperplanes is among the results of the paper \cite{GRV}, while at the present the preceding question seems to us a difficult problem.

\section*{Acknowledgements}
This paper is based up a series of lectures the author gave at Washington University in St.~Louis, the University of Michigan, the University of Napoli Federico II, Purdue University, Notre Dame University, Indiana University, St.~Louis University, Rutgers University, and Georgetown University while on sabbatical at Washington University in St. Louis.  It is a pleasure to thank Steven Krantz and the mathematics department at Washington University in St. Louis for hosting him during the undertaking of this work, to Jos\'e Pastrana for his questions that gave further impetus for the present structure, to Cody Stockdale for discussions concerning Calder\'on and Zygmund's decomposition lemma, to Jean Van Schaftingen for comments on a draft of the manuscript, and to Bogdan Rai\cb{t}\v{a} for his clarifications concerning a number of points in Section \ref{openproblems}.  Needless to say that I remain responsible for the remaining shortcomings.  The author is supported in part by the Taiwan Ministry of Science and Technology under research grants 105-2115-M-009-004-MY2, 107-2918-I-009-003 and 107-2115-M-009-002-MY2.



\end{document}